# A GENERAL FAMILY OF ESTIMATORS FOR ESTIMATING POPULATION MEAN USING KNOWN VALUE OF SOME POPULATION PARAMETER(S)

**Dr. M. Khoshnevisan**
GBS, Giffith University, Australia-(m.khoshnevisan@griffith.edu.au)

**Dr. Rajesh Singh, Pankaj Chauhan, Nirmala Sawan**
School of Statistics, DAVV, Indore (M.P.), India (rsingh.stat@dauniv.ac.in )

**Dr. Florentin Smarandache**
University of New Mexico, USA (smarand@unm.edu)

**Abstract :** A general family of estimators for estimating the population mean of the variable under study, which make use of known value of certain population parameter(s), is proposed. Under Simple Random Sampling Without Replacement (SRSWOR) scheme, the expressions of bias and mean-squared error (MSE) up to first order of approximation are derived. Some well known estimators have been shown as particular member of this family. An empirical study is carried out to illustrate the performance of the constructed estimator over others.

**Keywords** : Auxiliary information, general family of estimators, bias, mean-squared error, population parameter(s).

1. **Introduction**

Let y and x be the real valued functions defined on a finite population $U = (U_1, U_2, \ldots, U_N)$ and $\bar{Y}$ and $\bar{X}$ be the population means of the study character y and auxiliary character x respectively. Consider a simple random sample of size n drawn without replacement from population U. In order to have a survey estimate of the population mean $\bar{Y}$ of the study character y, assuming the knowledge of population mean $\bar{X}$ of the auxiliary character x, the well-known ratio estimator is

$$t_1 = \bar{y}\frac{\bar{X}}{\bar{x}} \tag{1.1}$$

Product method of estimation is well-known technique for estimating the populations mean of a study character when population mean of an auxiliary character is known and it is negatively correlated with study character. The conventional product estimator for $\bar{Y}$ is defined as

$$t_2 = \bar{y}\frac{\bar{x}}{\bar{X}} \tag{1.2}$$

Several authors have used prior value of certain population parameters (s) to find more precise estimates. Searls (1964) used Coefficient of Variation (CV) of study character at estimation stage. In practice this CV is seldom known. Motivated by Searls (1964) work, Sisodiya and Dwivedi (1981) used the known CV of the auxiliary character for estimating population mean of a study character in ratio method of estimation. The use of prior value of Coefficient of Kurtosis in estimating the population variance of study character y was first made by Singh et.al.(1973). Later, used by Sen (1978), Upadhyaya and Singh (1984) and Searls and Interpanich (1990) in the estimation of population mean of study character. Recently Singh and Tailor (2003) proposed a modified ratio estimator by using the known value of correlation coefficient.

In this paper, under SRSWOR, we have suggested a general family of estimators for estimating the population mean $\bar{Y}$. The expressions of bias and MSE, up to the first order of approximation, have been obtained, which will enable us to obtain the said expressions for any member of this family. Some well known estimators have been shown as particular member of this family.

## 2. The suggested family of estimators-

Following Walsh (1970), Reddy (1973) and Srivastava (1967), we define a family of estimators $\bar{Y}$ as

$$t = \bar{y}\left[\frac{a\bar{X}+b}{\alpha(a\bar{x}+b)+(1-\alpha)(a\bar{X}+b)}\right]^g \tag{2.1}$$

where a($\neq 0$), b are either real numbers or the functions of the known parameters of the auxiliary variable x such as standard deviation ($\sigma_x$), Coefficients of Variation ($C_X$), Skewness ($\beta_1(x)$), Kurtosis ($\beta_2(x)$) and correlation coefficient ($\rho$).

To obtain the bias and MSE of t, we write

$$\bar{y} = \bar{Y}(1+e_0), \bar{x} = \bar{X}(1+e_1)$$

such that

E ($e_0$)=E ($e_1$)=0,

and

$$E(e_0^2) = f_1 C_y^2, E(e_1^2) = f_1 C_x^2, E(e_0 e_1) = f_1 \rho C_y C_x,$$

where

$$f_1 = \frac{N-n}{nN}, \quad C_y^2 = \frac{S_y^2}{\bar{Y}^2}, \quad C_x^2 = \frac{S_x^2}{\bar{X}^2}.$$

Expressing t in terms of e's, we have

$$t = \bar{Y}(1+e_0)(1+\alpha\lambda e_1)^{-g} \tag{2.2}$$

where $\lambda = \dfrac{a\bar{X}}{a\bar{X}+b}$. (2.3)

We assume that $|\alpha\lambda e_1| < 1$ so that $(1+\alpha\lambda e_1)^{-g}$ is expandable.

Expanding the right hand side of (2.2) and retaining terms up to the second powers of e's, we have

$$t = \bar{Y}\left[1 + e_0 - \alpha\lambda g e_1 + \dfrac{g(g+1)}{2}\alpha^2\lambda^2 e_1^2 - \alpha\lambda g e_0 e_1\right] \quad (2.4)$$

Taking expectation of both sides in (2.4) and then subtracting $\bar{Y}$ from both sides, we get the bias of the estimator t, up to the first order of approximation, as

$$B(t) = f_1\bar{Y}\left[\dfrac{g(g+1)}{2}\alpha^2\lambda^2 C_x^2 - \alpha\lambda g \rho C_y C_x\right] \quad (2.5)$$

From (2.4), we have

$$(t - \bar{Y}) \cong \bar{Y}[e_0 - \alpha\lambda g e_1] \quad (2.6)$$

Squaring both sides of (2.6) and then taking expectations, we get the MSE of the estimator t, up to the first order of approximation, as

$$MSE(t) = f_1\bar{Y}^2\left[C_y^2 + \alpha^2\lambda^2 g^2 C_x^2 - 2\alpha\lambda g \rho C_y C_x\right] \quad (2.7)$$

Minimization of (2.7) with respect to $\alpha$ yields its optimum value as

$$\alpha = \dfrac{K}{\lambda g} = \alpha_{opt} \text{ (say)} \quad (2.8)$$

where

$$K = \rho\dfrac{C_y}{C_x}.$$

Substitution of (2.8) in (2.7) yields the minimum value of MSE (t) as

$$\min.MSE(t) = f_1\bar{Y}^2 C_y^2(1-\rho^2) = MSE(t)_0 \quad (2.9)$$

The min. MSE (t) at (2.9) MSE (t) is same as that of the approximate variance of the usual linear regression estimator.

## 3. Some members of the proposed family of the estimators' t

The following scheme presents some of the important known estimators of the population mean which can be obtained by suitable choice of constants $\alpha$, a and b:

| Estimator | Values of | | | |
|---|---|---|---|---|
| | $\alpha$ | a | b | g |
| 1. $t_0 = \bar{y}$<br>The mean per unit estimator | 0 | 0 | 0 | 0 |
| 2. $t_1 = \bar{y}\left(\dfrac{\bar{X}}{\bar{x}}\right)$<br>The usual ratio estimator | 1 | 1 | 0 | 1 |
| 3. $t_2 = \bar{y}\left(\dfrac{\bar{x}}{\bar{X}}\right)$<br>The usual product estimator | 1 | 1 | 0 | -1 |
| 4. $t_3 = \bar{y}\left(\dfrac{\bar{X}+C_x}{\bar{x}+C_x}\right)$<br>Sisodia and Dwivedi (1981) estimator | 1 | 1 | $C_x$ | 1 |
| 5. $t_4 = \bar{y}\left(\dfrac{\bar{x}+C_x}{\bar{X}+C_x}\right)$<br>Pandey and Dubey (1988) estimator | 1 | 1 | $C_x$ | -1 |

| Estimator | | | | | |
|---|---|---|---|---|---|
| 6. $t_5 = \bar{y}\left[\dfrac{\beta_2(x)\bar{x} + C_x}{\beta_2(x)\bar{X} + C_x}\right]$<br><br>Upadhyaya and Singh (1999) estimator | 1 | $\beta_2(x)$ | $C_x$ | -1 | |
| 7. $t_6 = \bar{y}\left[\dfrac{C_x\bar{x} + \beta_2(x)}{C_x\bar{X} + \beta_2(x)}\right]$<br><br>Upadhyaya, Singh (1999) estimator | 1 | $C_x$ | $\beta_2(x)$ | -1 | |
| 8. $t_7 = \bar{y}\left[\dfrac{\bar{x} + \sigma_x}{\bar{X} + \sigma_x}\right]$<br><br>G.N.Singh (2003) estimator | 1 | 1 | $\sigma_x$ | -1 | |
| 9. $t_8 = \bar{y}\left[\dfrac{\beta_1(x)\bar{x} + \sigma_x}{\beta_1(x)\bar{X} + \sigma_x}\right]$<br><br>G.N.Singh (2003) estimator | 1 | $\beta_1(x)$ | $\sigma_x$ | -1 | |
| 10. $t_9 = \bar{y}\left[\dfrac{\beta_2(x)\bar{x} + \sigma_x}{\beta_2(x)\bar{X} + \sigma_x}\right]$<br><br>G.N.Singh (2003) estimator | 1 | $\beta_2(x)$ | $\sigma_x$ | -1 | |
| 11. $t_{10} = \bar{y}\left[\dfrac{\bar{X} + \rho}{\bar{x} + \rho}\right]$<br><br>Singh, Tailor (2003) estimator | 1 | 1 | $\rho$ | 1 | |
| 12. $t_{11} = \bar{y}\left[\dfrac{\bar{x} + \rho}{\bar{X} + \rho}\right]$<br><br>Singh, Tailor (2003) estimator | 1 | 1 | $\rho$ | -1 | |

| | | | | |
|---|---|---|---|---|
| 13. $t_{12} = \bar{y}\left[\dfrac{\bar{X}+\beta_2(x)}{\bar{x}+\beta_2(x)}\right]$ Singh et.al. (2004) estimator | 1 | 1 | $\beta_2(x)$ | 1 |
| 14. $t_{13} = \bar{y}\left[\dfrac{\bar{x}+\beta_2(x)}{\bar{X}+\beta_2(x)}\right]$ Singh et.al. (2004) estimator | 1 | 1 | $\beta_2(x)$ | -1 |

In addition to these estimators a large number of estimators can also be generated from the proposed family of estimators t at (2.1) just by putting values of $\alpha$, g, a, and b.

It is observed that the expression of the first order approximation of bias and MSE/Variance of the given member of the family can be obtained by mere substituting the values of $\alpha$, g, a and b in (2.5) and (2.7) respectively.

## 4. Efficiency Comparisons

Up to the first order of approximation, the variance/MSE expressions of various estimators are:

$$V(t_0) = f_1 \bar{Y}^2 C_y^2 \tag{4.1}$$

$$MSE(t_1) = f_1 \bar{Y}^2 \left[C_y^2 + C_x^2 - 2\rho C_y C_x\right] \tag{4.2}$$

$$MSE(t_2) = f_1 \bar{Y}^2 \left[C_y^2 + C_x^2 + 2\rho C_y C_x\right] \tag{4.3}$$

$$MSE(t_3) = f_1 \bar{Y}^2 \left[C_y^2 + \theta_1^2 C_x^2 - 2\theta_1 \rho C_y C_x\right] \tag{4.4}$$

$$MSE(t_4) = f_1 \bar{Y}^2 \left[C_y^2 + \theta_1^2 C_x^2 + 2\theta_1 \rho C_y C_x\right] \tag{4.5}$$

$$MSE(t_5) = f_1 \bar{Y}^2 \left[C_y^2 + \theta_2^2 C_x^2 + 2\theta_2 \rho C_y C_x\right] \tag{4.6}$$

$$MSE(t_6) = f_1 \bar{Y}^2 \left[C_y^2 + \theta_3^2 C_x^2 + 2\theta_3 \rho C_y C_x\right] \tag{4.7}$$

$$MSE(t_7) = f_1\bar{Y}^2\left[C_y^2 + \theta_4^2 C_x^2 + 2\theta_4 \rho C_y C_x\right] \quad (4.8)$$

$$MSE(t_8) = f_1\bar{Y}^2\left[C_y^2 + \theta_5^2 C_x^2 + 2\theta_5 \rho C_y C_x\right] \quad (4.9)$$

$$MSE(t_9) = f_1\bar{Y}^2\left[C_y^2 + \theta_6^2 C_x^2 + 2\theta_6 \rho C_y C_x\right] \quad (4.10)$$

$$MSE(t_{10}) = f_1\bar{Y}^2\left[C_y^2 + \theta_7^2 C_x^2 - 2\theta_7 \rho C_y C_x\right] \quad (4.11)$$

$$MSE(t_{11}) = f_1\bar{Y}^2\left[C_y^2 + \theta_7^2 C_x^2 + 2\theta_7 \rho C_y C_x\right] \quad (4.12)$$

$$MSE(t_{12}) = f_1\bar{Y}^2\left[C_y^2 + \theta_8^2 C_x^2 - 2\theta_8 \rho C_y C_x\right] \quad (4.13)$$

$$MSE(t_{13}) = f_1\bar{Y}^2\left[C_y^2 + \theta_8^2 C_x^2 + 2\theta_8 \rho C_y C_x\right] \quad (4.14)$$

where

$$\theta_1 = \frac{\bar{X}}{\bar{X}+C_x}, \theta_2 = \frac{\beta_2(x)\bar{X}}{\beta_2(x)+C_x}, \theta_3 = \frac{C_x \bar{X}}{C_x \bar{X}+C_x},$$

$$\theta_4 = \frac{\bar{X}}{\bar{X}+\sigma_x}, \theta_5 = \frac{\beta_1(x)\bar{X}}{\beta_1(x)+\sigma_x}, \theta_6 = \frac{\beta_2(x)\bar{X}}{\beta_2(x)\bar{X}+\sigma_x},$$

$$\theta_7 = \frac{\bar{X}}{\bar{X}+\rho}, \theta_8 = \frac{\bar{X}}{\bar{X}+\beta_2(x)}.$$

To compare the efficiency of the proposed estimator t with the existing estimators $t_0$-$t_{13}$, using (2.9) and (4.1)-(4.14), we can, after some algebra, obtain

$$V(t_0) - MSE(t)_0 = C_y^2 \rho^2 > 0 \quad (4.15)$$

$$MSE(t_1) - MSE(t)_0 = (C_x - \rho C_y)^2 > 0 \quad (4.16)$$

$$MSE(t_2) - MSE(t)_0 = (C_x + \rho C_y)^2 > 0 \quad (4.17)$$

$$MSE(t_3) - MSE(t)_0 = (\theta_1 C_x - \rho C_y)^2 > 0 \quad (4.18)$$

$$MSE(t_4) - MSE(t)_0 = (\theta_1 C_x + \rho C_y)^2 > 0 \quad (4.19)$$

$$MSE(t_5) - MSE(t)_0 = (\theta_2 C_x + \rho C_y)^2 > 0 \qquad (4.20)$$

$$MSE(t_6) - MSE(t)_0 = (\theta_3 C_x + \rho C_y)^2 > 0 \qquad (4.21)$$

$$MSE(t_7) - MSE(t)_0 = (\theta_4 C_x + \rho C_y)^2 > 0 \qquad (4.22)$$

$$MSE(t_8) - MSE(t)_0 = (\theta_5 C_x + \rho C_y)^2 > 0 \qquad (4.23)$$

$$MSE(t_9) - MSE(t)_0 = (\theta_6 C_x + \rho C_y)^2 > 0 \qquad (4.24)$$

$$MSE(t_{10}) - MSE(t)_0 = (\theta_7 C_x - \rho C_y)^2 > 0 \qquad (4.25)$$

$$MSE(t_{11}) - MSE(t)_0 = (\theta_7 C_x + \rho C_y)^2 > 0 \qquad (4.26)$$

$$MSE(t_{12}) - MSE(t)_0 = (\theta_8 C_x - \rho C_y)^2 > 0 \qquad (4.27)$$

$$MSE(t_{13}) - MSE(t)_0 = (\theta_8 C_x + \rho C_y)^2 > 0 \qquad (4.28)$$

Thus from (4.15) to (4.28), it follows that the proposed family of estimators 't' is more efficient than other existing estimators $t_0$ to $t_{13}$. Hence, we conclude that the proposed family of estimators 't' is the best (in the sense of having minimum MSE).

## 5. Numerical illustrations

We consider the data used by Pandey and Dubey (1988) to demonstrate what we have discussed earlier. The population constants are as follows:

N=20, n=8, $\bar{Y} = 19.55$, $\bar{X} = 18.8$, $C_x^2 = 0.1555$, $C_y^2 = 0.1262$, $\rho_{yx} = -0.9199$, $\beta_1(x) = 0.5473$, $\beta_2(x) = 3.0613$, $\theta_4 = 0.7172$.

We have computed the percent relative efficiency (PRE) of different estimators of $\bar{Y}$ with respect to usual unbiased estimator $\bar{y}$ and compiled in table 5.1.

**Table 5.1: Percent relative efficiency of different estimators of $\bar{Y}$ with respect to $\bar{y}$**

| Estimator | PRE |
|---|---|
| $\bar{y}$ | 100 |
| $t_1$ | 23.39 |
| $t_2$ | 526.45 |
| $t_3$ | 23.91 |
| $t_4$ | 550.05 |
| $t_5$ | 534.49 |
| $t_6$ | 582.17 |
| $t_7$ | 591.37 |
| $t_8$ | 436.19 |
| $t_9$ | 633.64 |
| $t_{10}$ | 22.17 |
| $t_{11}$ | 465.25 |
| $t_{12}$ | 27.21 |
| $t_{13}$ | 644.17 |
| $t_{(opt)}$ | 650.26 |

From table 5.1, we observe that the proposed general family of estimators is preferable over all the considered estimators under optimum condition.